\documentclass{article}

\usepackage{amssymb}

\usepackage[T1]{fontenc}
\usepackage[dvips]{graphicx}

\newtheorem{theorem}{Theorem}[section]
\newtheorem{definition}[theorem]{Definition}
\newtheorem{lemma}[theorem]{Lemma}

\newcommand{\PROOF}{\noindent {\bf Proof}: }

\begin{document}

\title{Manifolds without $\frac{1}{k}$-geodesic}
\author{Wing Kai, Ho}
\date{\today}

\maketitle

\begin{abstract}
We will show that for any positive integer $\emph k$, there exist a
smooth manifold that has no $\frac{1}{\emph k}$-geodesic.
\end{abstract}

\section {Introduction}
For compact length spaces, it is known that the classical marked
length spectrum may not be continuous with respect to
Gromov-Hausdorff limit: there exists a sequence of manifolds $M_i$,
$M_i \rightarrow M$ in the Gromov-Hausdorff sense, such that closed
geodesics do not persist under this limit. C. Sormani has introduced
the $\frac {1}{k}~length~spectrum$ $L_{\frac{1}{k}}(M)$, the set of
lengths of $\frac{1}{k}-geodesics$ in $M$. A closed geodesic of
length $l$ is called a $\frac {1}{k}-geodesic$ if it is length
minimizing on every segment of length $\frac {l}{k}$. Sormani proved
that $\frac{1}{k}-geodesics$ persist under Gromov-Hausdorff limit,
which implies that $\frac {1}{k}~length~spectra$ are stable under
Gromov-Hausdorff convergence. For discussions about $\frac
{1}{k}~length~spectrum$,
see \cite{So}.\\

Sormani showed that every shortest homotopically non-trivial closed
geodesic is a $\frac{1}{k}-geodesic$, $\forall k \in \mathbb N$ [see
example below]. Sormani then proposed the following question: does
that exist $k \in \mathbb N$, such that every smooth, compact,
simply connected manifold has a $\frac{1}{k}-geodesic$? We address
this question by constructing metrics $\rho_k$ on $S^2$ for each $k
\in \mathbb N$, such that $(S^2,\rho_k)$ has no
$\frac{1}{k}-geodesic$. However it is still not known if every
manifold admits a metric without $\frac {1}{2}-geodesic$.

\section{Definition}

Let $M$ be a smooth manifold. $\gamma$: $S^1 \rightarrow$ $M$ be a
closed geodesic parameterized by arc length and have length $l$.
$\gamma$ is called $\frac{1}{2}-geodesic$ if it is distance
minimizing on every segment of length $\frac {\emph l}{2}$.
Similarly, a $\frac{1}{\emph k}-geodesic$ is a closed geodesic that
is distance minimizing on every segment of length $\frac {l}{k}$.

\vspace{0.2in}

$\bf Example :$ Suppose that $M$ is not simply connected. Let
$\gamma$ be a shortest homotopically non-trivial closed curve in
$M$. Then $\gamma$ is a closed geodesic (for instance, see
\cite{Jo}). Let us show that $\gamma$ is a $\frac{1}{2}$-geodesic.
Denote the length of $\gamma$ by $l$. Reasoning by contradiction,
assume that there are two points $p$, $q$ on $\gamma$ that are
$\frac{l}{2}$ apart along $\gamma$ and that can be connected by a
geodesic $\gamma_1$ that is shorter than $\frac{l}{2}$. The points
$p$ and $q$ divide $\gamma$ into two geodesics. Each of them can be
closed up by adding $\gamma_1$. Hence we represented $\gamma$ as a
product of two loops, each of which is shorter than $l$. Since
$\gamma$ is homotopically non-trivial, so is at least one of these
loops. This contradicts to our assumption that $\gamma$ is a
shortest homotopically non-trivial loop.\\

By a $segment$ of a geodesic $\gamma$ we mean the restriction of
$\gamma$ to a closed interval. A $\frac{1}{k}~segment$ is a segment
of length $\frac{l}{k}$. A $loop$ is a finite union of segments that
bound a 2 dimensional disc.

\section{Construction of the surfaces}

Our goal is to show that, for every integer $k \geqslant 2$, there
exists a smooth surface $M_k$ that has no $\frac {1}{k}$-geodesic.
In our construction, each $M_k$ will be a surface of revolution. We
first start with $k=2$, and then generalize to all $k$.

\vspace{0.2in}

{\bf The surface.} Consider a curve in ($\mathbb R^2$, Euclidean
metric) that consists of a straight line joining (0,1) and (n,0)
($n$ to be determined later), and a straight line from (0,1) to
(0,0). These are just two sides of a right triangle. If we revolve
this curve about the x-axis, we get a cone $K$ with circular base of
radius 1 and height $n$. Now smoothen the two angles on both ends of
the 'hypotenuse' by replacing a small neighborhood of each of the
angle with a smooth arc, so that when we revolve it about the x-axis
we get a smooth surface. The resulting surface is our $M_2$. For the
sake of simplicity, we create $M_2$ in the way that the longest
parallel (the great parallel) has radius 1. Now, $M_2$ is
diffeomorphic to $S^2$, and looks like a smoothened cone. Actually,
since we alter arbitrarily small neighborhoods of the angles, the
surface is 'metrically' close to $K$. i.e. there exists a map $f$
between $M_2$ and $K$, and $f$ has very small distortion. For
instance, such $f$ can be obtained by starting from the midpoint of
the
hypotenuse. We elongate it by sliding the two ends to sharp angles, followed by a suitable rescaling.\\

%figure one
\begin{figure}[h]
\caption{Construction of $M_k$}
\centerline{\resizebox{80mm}{!}{\includegraphics{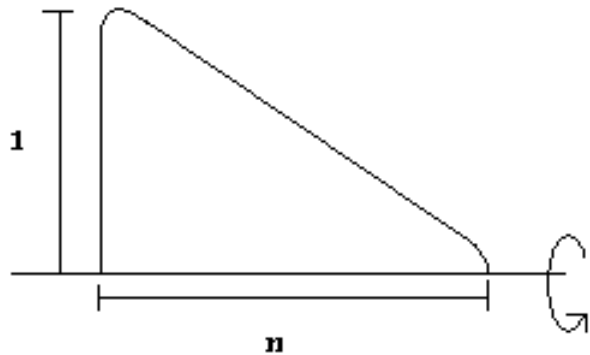}}}
\label{pic1}
\end{figure}

The rest of this section is dedicated to proving the following
statement:

\begin{theorem}
With $n$ suitably large, $M_2$ has no $\frac{1}{2}-geodesic$.
\end{theorem}
%insert figure here

\vspace{0.3in}

To prove our claim, we will show that all closed geodesics in $M_2$
are not $\frac{1}{2}-geodesic$. We begin with the following
observation:

\begin{lemma}
$\frac {1}{2}$-geodesic has no self-intersection.
\end{lemma}

\PROOF Suppose a closed geodesic $\gamma$ of length $l$ has
self-intersection. Then there exists a segment $\eta$ with two
endpoints coincide, such that $\eta$ has length $\leqslant
\frac{l}{2}$. To see this, suppose $\gamma$ has at least one
self-intersection. Then this self-intersection separates $\gamma$
into two geodesics, such that the four endpoints coincide at one
point. (Think of the figure 8). It is easy to see that one of them
has to have length less or equal to $\frac{l}{2}$. Now, any segment
of length $\frac{l}{2}$ that contains $\eta$ cannot be distance
minimizing. That is because the two endpoints of this segment can be
joined by a shorter path, obtained by deleting $\eta$ from the
segment.
$\square$\\

The reason that we consider surfaces of revolution is we can
classify all geodesics using $Clairaut's~integral$: Given a
geodesic, if we denote by $r$ the radius of the parallel which the
geodesic intersects with,
$\theta$ be the angle of intersection. Then the relation\\

\begin{equation}
r~\cos \theta = const =c,
\end{equation}

holds on the whole geodesic.\\

With this we have the following lemma:

\begin{lemma}
No closed geodesic can stay on one side of the great parallel (the
longest parallel). i.e. it must intersect the great parallel.
\end{lemma}
\PROOF  Arguing by contradiction, suppose the geodesic $\gamma$
stays on one side. By compactness of $\gamma$, there exist a
shortest and longest parallel (with radius $r_1$ and $r_2$), such
that $\gamma$ is tangential to both and lies between them. If
$r_1=r_2$, then $\gamma$ is a parallel. This cannot happen, since
any parallel of this kind is generated by the rotation of a point of
the profile curve where the tangent is not parallel to the axis of
revolution. None of these parallel can be geodesic \cite{Do}.
Therefore we must have $r_1 \neq r_2$. This contradicts the
$Clairaut's~integral$ since in this case, $c=r_1$ and $c=r_2$.
 $\square$\\

So any geodesic is uniquely determined by the following data: the
point of intersection with the great parallel and the angle of
intersection $\alpha$. Now by Clairaut's integral, the angle
$\alpha$ determines the constant $c=c_\alpha$. Denote this geodesic
$\gamma_{\alpha}(t)$, $\gamma_{\alpha}(0)$=point of intersection
with the great parallel.\\

Let's investigate all closed geodesics in $M_2$:\\
\vspace{0.1in}\\
$\bf {Meridians~(\alpha=\frac{\pi}{2})}:$ Meridians cannot be $\frac
{1}{2}-geodesic$ if $n$ is large enough. To see this, fix any
meridian, pick two points $p,q$ that lie on the same parallel and
split the meridian into halves. The distance between $p$ and $q$ is
approximately half of the
length of the parallel and thus is much shorter than the length of half-meridian.\\

\noindent $\bf Great~parallel~(\alpha=0):$ The longest parallel
(with radius 1) of $M_2$ cannot be $\frac {1}{2}-geodesic$. Fix any
two antipodal point $p,q$ on the great parallel. The distance
between $p$ and $q$ along the parallel is $\pi$. However $p$ and $q$
can be joined by a path across the base. The length of this path
equals approximately the diameter of the great parallel. Which means
$p$ and $q$ can be joined by a shorter path. Hence the great
parallel is not a
$\frac{1}{2}-geodesic$.\\

\noindent $\bf Other~closed~geodesics~(\alpha \in
(0,\frac{\pi}{2})):$ Without loss of generality, we can assume
$\gamma_\alpha'(0)$ is pointing into the cone. Let $r_{\alpha}(\emph
t)$ be the radius of parallel intersecting $\gamma_{\alpha}$ at
$\gamma_{\alpha}($\emph t$)$, and $\theta_{\alpha}(\emph t)$ be
angle of intersection. Observe that when
$r_{\alpha}(t_\alpha)=c_{\alpha}$, for some $t_\alpha \in [0,l]$,
$\gamma_{\alpha}$ is tangential to the parallel, and then it will
start to return~\cite{Sp}. Denote by $R_{\alpha}$ the parallel where
$\gamma_{\alpha}$ start to turn back.\\

\begin{definition}
For each $\alpha \in [0,\frac{\pi}{2})$, define the total rotation
$T_\alpha(t),~t \in [0,l]$ to be the net (oriented) angle of
rotation of $\gamma_\alpha$ about the axis of revolution from
$\gamma_\alpha(0)$ to $\gamma_\alpha(t)$.
\end{definition}

Example: When $\alpha$=0, $\gamma_\alpha$ is just the great
parallel, Therefore $T_\alpha(t)=\pm t$ (depending on the orientation chosen)\\

\vspace{0.1in} Firstly, for any $\alpha \neq \frac{\pi}{2}$,
$|T_\alpha(t)|$ is a monotonic increasing function. This is
equivalent to saying that any non-meridian geodesic $\gamma$ rotates
only in one direction. To prove this claim, assume on the contrary
that $\gamma$ changes rotational direction at some point. Then at
this point, $\gamma$ should be tangential to a meridian. By the
uniqueness of geodesics (in a smooth manifold, a point and a vector
uniquely determine a geodesic), $\gamma$ should coincide with a
meridian. This contradicts the
assumption that $\gamma$ is a non-meridian.\\

Now recall that $\gamma_\alpha (t_\alpha)$ is the point when
$\gamma_\alpha$ turns back, we have the following lemma:

\begin{lemma}
If $|T_\alpha(t_\alpha)|>\pi$, then $\gamma_\alpha$ has
self-intersection.
\end{lemma}
\PROOF As noted above, $\gamma_\alpha$ only rotates in one
direction, so from the same analysis on the $Clairaut's~integral$,
when $\gamma_\alpha$ returns and hits the great parallel again, it
should have rotated by $|2T_\alpha(t_\alpha)|$. We know from
$Clairaut's~integral$ that $\gamma_\alpha$ cannot touch the great
parallel. Therefore, if we let $x$ to be the center of the base, the
winding number of $\gamma_\alpha$ about $x$ is strictly greater than
1 (up to orientation). Which implies that $\gamma_\alpha$ has
self-intersection. $\square$\\

We are now ready to list all the remaining geodesics in $M_2$, to
simplify our argument, let us divide $M_2$ into four areas. Recall
that in our construction, we smoothen 2 corners of the generating
curve. Therefore when we revolve it: There is a curved cap at the
tip (the cap), a thin curved belt around the great parallel (the
belt), a flat disc at the bottom (the disc) and the long cone (the
cone) [figure 2]. Only the cap and the belt have non-zero curvature.\\

The remaining geodesics can be divided into
three types:\\
a) Geodesics that never leave the belt.\\
b) Geodesics that enter the cap.\\
c) Geodesics that enter the cone but miss the cap. \\
\vspace{0.1in}\\

%figure
\begin{figure}[h]
\caption{Four areas of $M_2$}
\centerline{\resizebox{80mm}{!}{\includegraphics{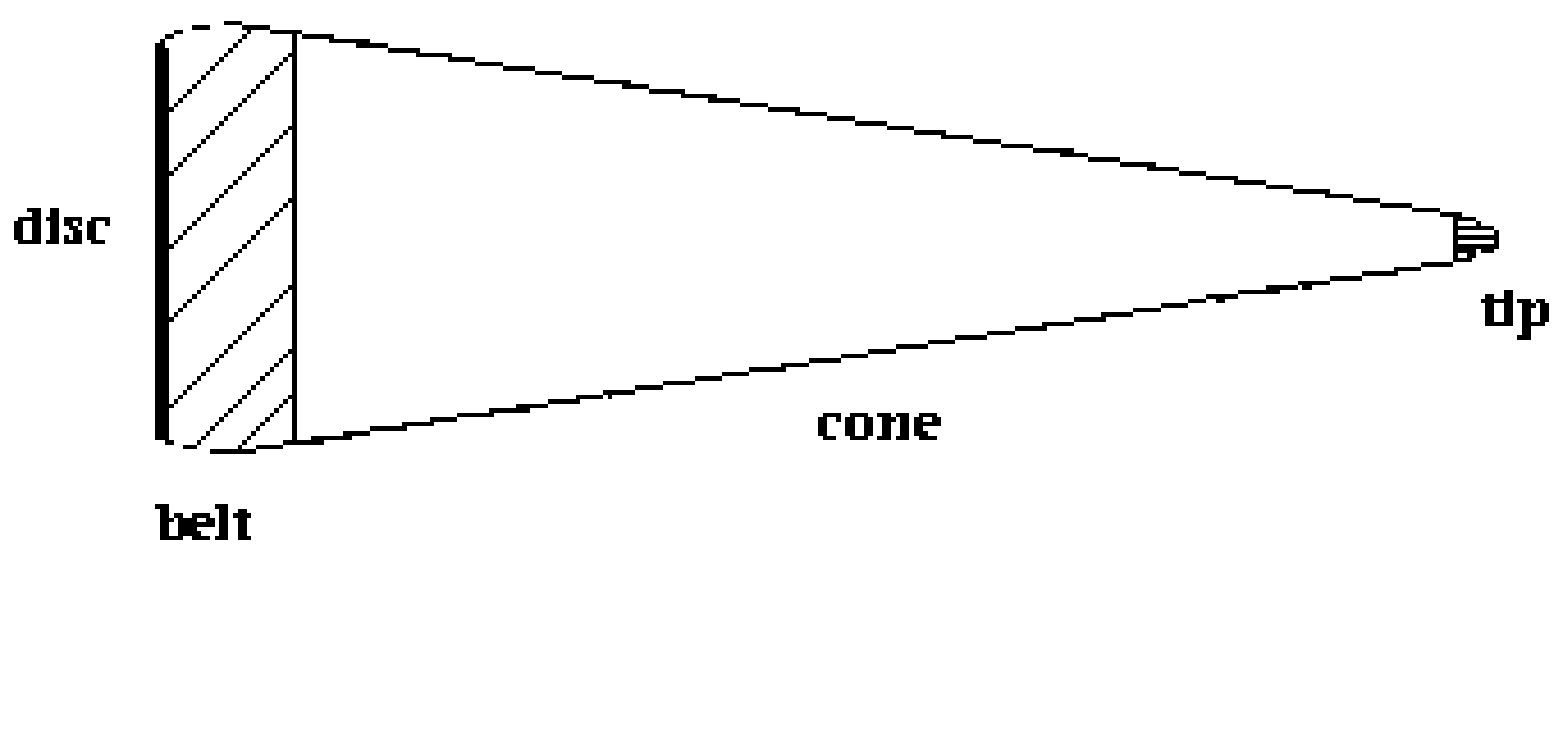}}}
\label{pic2}
\end{figure}

There are two parallels which separate the belt and the cone, the
cone and the cap. Denote these two parallels by $R_{\alpha'}$ and
$R_{\alpha''}$ respectively. Recall that $R_{\alpha}$ is the
parallel where $\gamma_\alpha$ starts to turn back. Now since in
constructing $M_2$, the belt and the cap can be arbitrarily thin. We
can assume $\alpha', (\frac {\pi}{2}-\alpha'') \ll \frac{\pi}{2}$.
To make the following arguments simpler, we also dilate $M_2$
proportionally so that $R_{\alpha'}$ has length 1. This has no
impact to all previous arguments.\\

The three cases of geodesics are equivalent to:\\

\noindent a) $\alpha \in (0,\alpha')$\\
b) $\alpha \in [\alpha'',\frac {\pi}{2})$\\
c) $\alpha \in [\alpha',\alpha'')$\\

{\bf Case a)} If $\gamma_\alpha$ wraps around $M_2$ twice or more,
by the same argument of winding number, $\gamma_\alpha$ has
self-intersection. Hence by lemma 3.2, $\gamma_\alpha$ is not
$\frac{1}{2}-geodesic$. If $\gamma_\alpha$ only wrap around $M_2$
once, then it enters each side of the great parallel once. Suppose
the distance between $R_{\alpha'}$ and the great parallel is
$\epsilon$, $\epsilon \ll 1$. Then $\gamma_\alpha$'s length is
within $2\pi \pm 10\epsilon$. Therefore $\gamma_\alpha$ is similar
to the great parallel: any two points $p,q$ dividing $\gamma_\alpha$
into halves can be joined by a path of length $\leq 2+10\epsilon$.
This is a shorter path. Therefore we conclude that all geodesics in
this case are not
$\frac{1}{2}-geodesic$.\\

{\bf Case b)} Now, since $\gamma_\alpha$ connects the great parallel
and some point in the cap, $\gamma_\alpha$ has at least length of
($2n-\epsilon'$) for some small $\epsilon'$. Then it is just like a
meridian: find two points which are $\frac{2n-\epsilon'}{2}$ apart
and lie on the same parallel. When $n$ is large the half-parallel is
a shorter path. Hence no
geodesic in case b can be $\frac{1}{2}-geodesic$.\\

{\bf Case c)} If $\gamma_\alpha$ enters the cone, then it must cross
the parallel $R_{\alpha'}$. So there is an angle of intersection
$\tilde{\alpha}$ between $\gamma_\alpha$ and $R_{\alpha'}$. Define
$T(\tilde{\alpha})$, the $first~return~ rotation$ to be the total
rotation of $\gamma_\alpha$ from $R_{\alpha'}$ and the point when it
first hit $R_{\alpha'}$ again [figure 3].

%figure 3
%figure
\begin{figure}[h]
\caption{$T(\tilde{\alpha})=2\pi$}
\centerline{\resizebox{80mm}{!}{\includegraphics{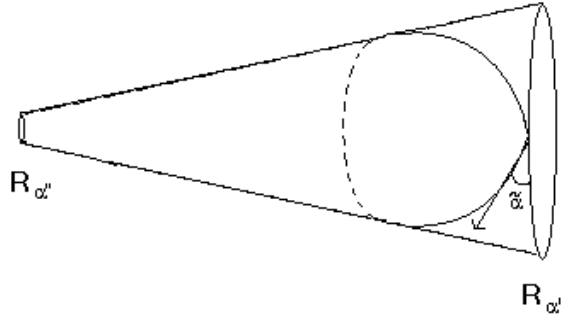}}}
\label{pic3}
\end{figure}

We need the following lemma:\\
\begin{lemma}
$T(\tilde{\alpha})$ is monotonic increasing in $\tilde{\alpha}$ for
all geodesics in case c.
\end{lemma}
\PROOF Consider the universal cover of the cone. Construct it by
starting with an annulus, cut through one radius. Then take another
copy of the same thing and then glue the left side of the cut from
the first copy to the right side of the second copy.  Continue
infinitely we get the universal cover. It looks like a infinite
spiral and is a topological infinite strip. A fundamental domain is a sector [figure 4].\\

%figure 4
\begin{figure}[h]
\caption{The universal cover of the cone}
\centerline{\resizebox{80mm}{!}{\includegraphics{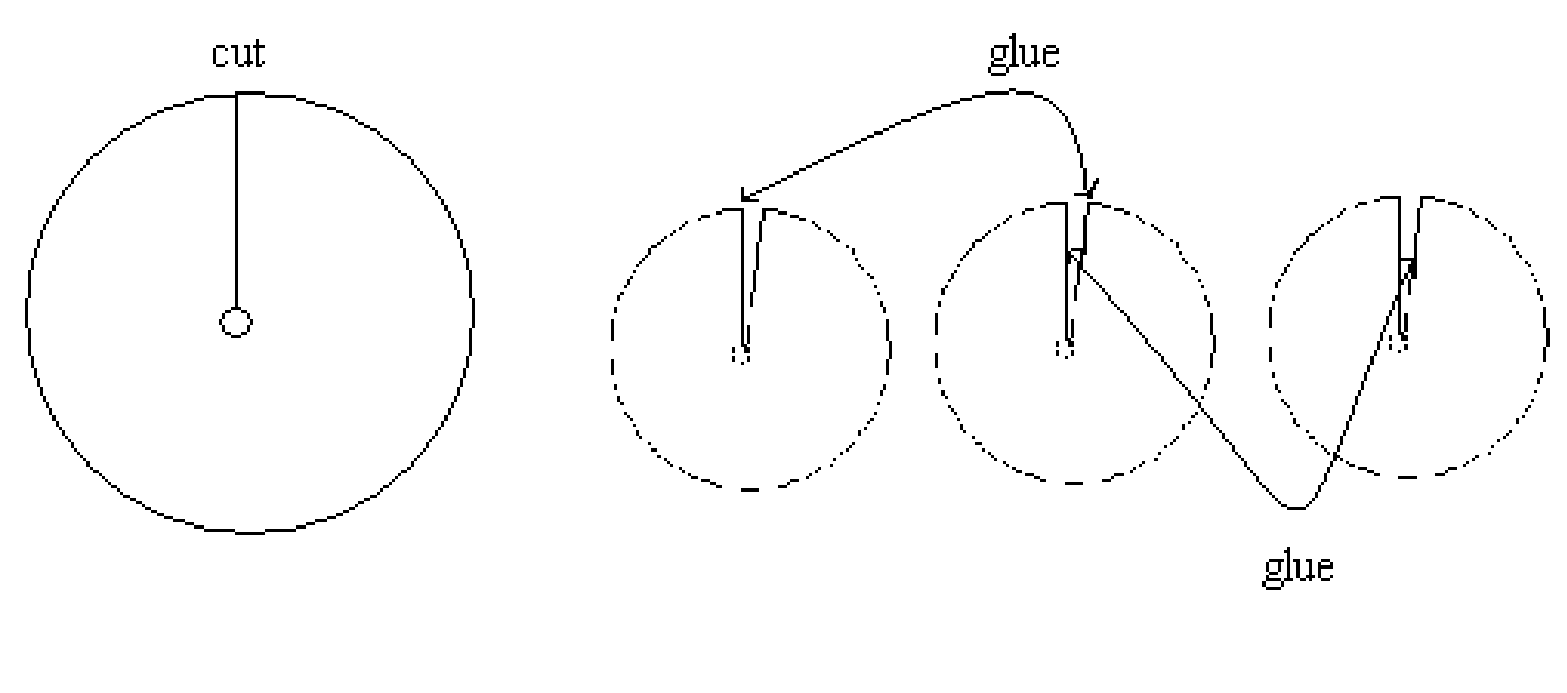}}}
\label{pic4}
\end{figure}

Now this is a development of the cone area, any geodesic segment is
a straight line. Also, $\tilde{\alpha}$ is given by the angle of
intersection with the outer circle. It is now easy to see that
$T(\tilde{\alpha})$ is monotonic increasing in $\tilde{\alpha}$:
Since we assume that $R_{\alpha'}$ has length 1, $T(\tilde{\alpha})$
is the length of the arc corresponding to the chord given by $\gamma_\alpha$ [figure 5]. $\square$\\

%figure 5
\begin{figure}[h]
\caption{$T(\tilde{\alpha})$ is monotonic increasing}
\centerline{\resizebox{80mm}{!}{\includegraphics{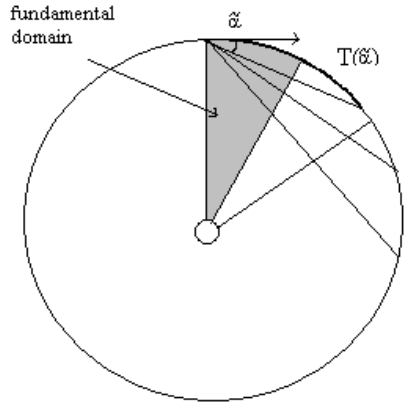}}}
\label{pic5}
\end{figure}

Finally, we claim that for any fix $\zeta \geq \epsilon$. When $n$
is large enough, any $\gamma_\alpha$ not contained in the
$\zeta$-neighborhood of the great parallel has self-intersection. To
see this, consider  the fundamental domain (with arc length
$l(R_{\alpha'})=1$). A chord connecting two end points of the arc is
a geodesic $\gamma_\alpha$ with $T(\tilde{\alpha})=2\pi$. Denote by
$L$ the distance between $R_{\alpha'}$ and $\gamma_\alpha$.
Elementary calculation shows that
$L=n(1-\sqrt{1-\sin^2\frac{1}{2n}}) \rightarrow 0$ as $n \rightarrow
\infty$ [figure 6]. So when $n$ is large enough such that $L=\zeta$,
the geodesic that turns back exactly at the boundary of the
$\zeta$-neighborhood gives $T(\tilde{\alpha})=2\pi$, hence it has
self-intersection. Together with Lemma 3.6., when $\gamma_\alpha$ is
not contained in the $\zeta$-neighborhood of the great parallel, it
has self-intersection. Therefore by Lemma 3.2, such
geodesic cannot be $\frac{1}{2}-geodesic$.\\

%figure 6
\begin{figure}[h]
\caption{$L \rightarrow 0$ as $n \rightarrow \infty$}
\centerline{\resizebox{80mm}{!}{\includegraphics{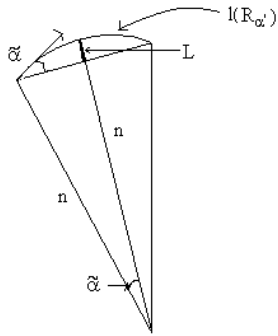}}}
\label{pic6}
\end{figure}

Now, the remaining geodesics are those that sit inside the
$\zeta$-neighborhood of the great parallel. Take $\zeta \ll 1$, this
is a similar case as the geodesics that is contained in the curved
belt: any two points $p,q$ dividing $\gamma_\alpha$ into halves can
be
joined by a shorter path through the disc.\\

So if we choose $n$ large enough such that all the previous criteria
are
met. Then $M_2$ has no $\frac{1}{2}-geodesic$ and we finish the proof of Theorem 3.1.$\square$\\

\section{When $k \geq 3$}
The construction of $M_k$ is similar to that of $M_2$, except that
we have to use larger $n$, thinner belt and smaller cap.
\begin{theorem}
For any fixed $k$, $M_k$ has no $\frac{1}{k}-geodesic$.
\end{theorem}

As what we have done before, we will exhibit all possible geodesics.
First off, any closed geodesic $\gamma$ must intersect the great
parallel (Lemma 3.3). So as before we can use the angle of
intersection $\alpha$ to
characterize the geodesics. In the following we still assume that $\gamma$ has length $l$.\\

$\bf {Meridians} : $ Meridians are not $\frac {1}{k}-geodesic$ if
$n$ is large enough. Again, find two points $p,q$ near the tip that
contain a $\frac {1}{k}$ segment and lie on the same parallel. $n$
being large implies $\frac{l}{k}$ is much larger than the length of
any parallel. Therefore there is a shorter
path joining $p, q$.\\

$\bf {Great~parallel} : $ The great parallel has length 2$\pi$. Any
two points $p,q$ that contains a $\frac{1}{k}$ segment
($\frac{2\pi}{k}$ long) of the great parallel can be joined by a
shorter path through the base. This is a
chord on the disc plus some small error.\\

$\bf {Other~geodesics} : $ Again, these geodesics can be categorized
into 3 types: stays in the belt, goes into the cap, and goes into the cone but not the cap.\\

1) In the belt: If the geodesic wraps around once, then it is
similar to the case of the great parallel: $p,q$ can be joined by a
shorter path that is close to a chord of the great parallel. If the
geodesic wraps around $m$ times, then for $p,q$ bounding a
$\frac{1}{k}$ segment, they are apart by approximately $\frac
{2m\pi}{k} >\frac{2\pi}{k}$. Again, $p,q$ can be joined by a shorter
path through the disc.\\

2) Into the cap: Similar to the case of $k=2$, any geodesic that
runs into the cap has length at least $2n-\epsilon'$ for some small
$\epsilon'$. We can find $p,q$ near the tip. Such that $p,q$ bound a
$\frac {1}{k}$ segment ($\frac {2n-\epsilon'}{k}$ long) of the
geodesic, and lie on the same parallel. Then $p,q$ can be joined by
a
path close to a half-parallel which is a shorter path.\\

3) Geodesics that run into the cone but miss the cap: Since $k \geq
3$, lemma 3.2 no longer applies here. However, we have the
following lemma:\\

\begin{lemma}
For any $\gamma_\alpha$ in case 3. If $\gamma_\alpha$ has $(k+1)$
self-intersections in the cone area. Then $\gamma_\alpha$ is not a
$\frac{1}{k}-geodesic$.
\end{lemma}

\PROOF Suppose $\gamma_\alpha$ has $(k+1)$ self-intersections in the
cone area. Recall that by $Clairaut's~integral$, any geodesic of
this form is symmetric about the meridian that contains the point
where the geodesic starts to turn back. The self-intersections split
$\gamma_\alpha$ into at least $(2k+1)$ segments. Let us label the
corresponding segments 1, 2, 2', etc. [figure 7]. Notice that
segment 1 forms a loop, segments 2 and 2' form another loop and so
on. There are altogether $k$ loops of this kind in the cone area.

%figure 7
\begin{figure}[h]
\caption{A geodesic in case 3}
\centerline{\resizebox{80mm}{!}{\includegraphics{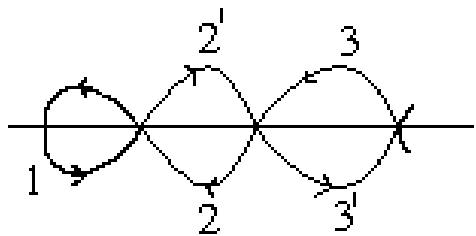}}}
\label{pic7}
\end{figure}

Now we consider the universal cover again. Since segment 1 is the
only one which is orthogonal to a meridian. It has to be strictly
shorter than $length(segment~i)+length(segment~i')$ for $2 \leq i
\leq k$ [figure 8]. That means segment 1 is the shortest loop among
the $k$ loops in the cone area. Which implies $length(segment~1)
<\frac{l}{k}$. Any $\frac{1}{k}$ segment of $\gamma_\alpha$
containing segment 1 cannot be shortest path. Since we can connect
the two endpoints by
a shorter path if we jump segment 1 at the point of intersection.$\square$\\

Note that in general $\frac{1}{n}-geodesic$ can have as many as $(\frac{n}{2}-1)^2$ self-intersections.\\

%figure 8
\begin{figure}[h]
\caption{Segment 1 has length $\leq \frac{l}{k}$}
\centerline{\resizebox{80mm}{!}{\includegraphics{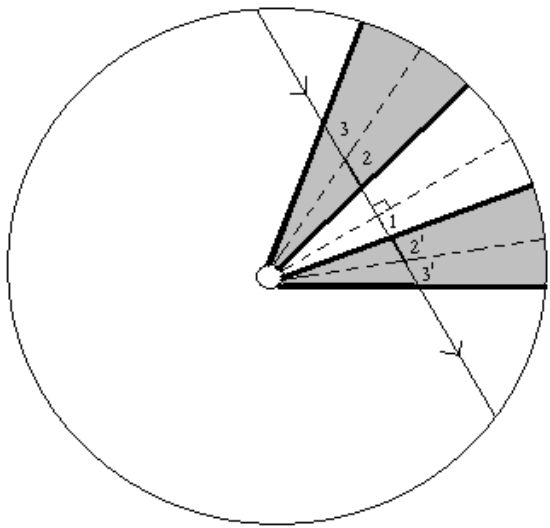}}}
\label{pic8}
\end{figure}

Now given any fixed $\zeta$, $\epsilon < \zeta \ll \frac{1}{k}$.
Using the same argument as $k=2$ [page 7]. When $n$ is large enough,
$T(\tilde{\alpha})>2(k+1)\pi$ for all $\gamma_\alpha$ not contained
inside the $\zeta$-neighborhood of the great parallel. This implies
that $\gamma_\alpha$ has $(k+1)$ self-intersections and by lemma
4.2, $\gamma_\alpha$ is not $\frac{1}{k}-geodesic.$ If
$\gamma_\alpha$ is contained inside the $\zeta$-neighborhood, then $
\zeta \ll \frac{1}{k}$ implies $\gamma_\alpha$ is similar to those
in case 1, hence it cannot be
$\frac{1}{k}-geodesic$.\\

So for $n$ large enough, $M_k$ has no $\frac{1}{k}-geodesic$.\\

To conclude, we have proven the following result:\\

\begin{theorem}
For any fixed $k$, there exist a metric $\rho_k$ on $S^2$ such that
$(S^2,\rho_k)$ has no $\frac{1}{k}-geodesic$.
\end{theorem}

\end{document}